\documentclass[leqno,12pt]{amsart}
\usepackage{amssymb} 

\theoremstyle{plain}
\newtheorem{theorem}{\indent\sc Theorem}[section]
\newtheorem{lemma}[theorem]{\indent\sc Lemma}
\newtheorem{corollary}[theorem]{\indent\sc Corollary}

\theoremstyle{definition}

\newtheorem*{remark0}{\indent\sc Remark}

\begin{document}
\title[Stability and Instability for Harmonic Foliations]{Foliations with complex leaves and instability for harmonic foliations}
\author[K. Ichikawa]{Kei Ichikawa}
\author[T. Noda]{Tomonori Noda} 
\address{Kawagoe Senior High School \endgraf
2-6 Kuruwa \endgraf
Kawagoe Saitama \endgraf
350-0053 \endgraf
Japan \endgraf
}

\address{Graduate school of science \endgraf
Osaka University \endgraf
1-1 Machikaneyama \endgraf
Toyonaka, Osaka \endgraf
560-0043 \endgraf
Japan}
\email{momentmap@yahoo.co.jp}
\maketitle

\begin{abstract}
In this paper, we study stability for harmonic foliations on locally conformal  K\"ahler manifolds with complex leaves. We also discuss instability for harmonic foliations on compact submanifolds immersed in Euclidean spaces and compact homogeneous spaces.
\end{abstract}

\section{Introduction} 

In this paper we study stability and instability for harmonic foliations. Let $(M,J,g_M)$ be a Hermitian manifold and $\Omega$ the fundamental $2$-form associated with $g_M$. Then $(M,J,g_M)$ is a {\it locally conformal K\"ahler} manifold if there exists a closed $1$-form $\omega$, called the {\it Lee form}, satisfying $d\Omega=\omega\wedge\Omega$. Besides K\"ahler manifolds, there are numerous examples of locally conformal K\"ahler manifolds. For instance, a Vaisman manifold is known to be a locally conformal K\"ahler manifold with non-exact and parallel Lee form.\vspace{4mm}

The main purpose of this paper is to prove the following stability theorem for harmonic foliations on compact locally conformal K\"ahler manifolds: \vspace{4mm}

{\sc Main Theorem.} {\it Let $(M,J,g_M)$ be an $n$-dimensional compact locally conformal K\"ahler manifold. If $\mathcal F$ is a harmonic foliation on $M$ with bundle-like metric $g_M$ foliated by complex submanifolds, then $\mathcal F$ is stable.}\vspace{1pc}

This is an analogoue of the theorem ``a holomorphic map between two K\"ahler manifolds is stable as a harmonic map'' (see also Corollary 3.1), where harmonicity for a foliation ${\mathcal F}$ on a Riemannian manifold $(N,g_N)$ is defined by Kamber and Tondeur in [6] as the harmonicity of the canonical projection $\pi$ from $TN$ onto the normal bundle $Q$ for the foliation ${\mathcal F}$. The key of the proof of Main Theorem is the compatibility of the complex structure with the connection on the normal bundle of the foliation (see Lemma 3.2).\vspace{4mm}

We also discuss instability for harmonic foliations on compact homogeneous Riemannian manifolds or compact submanifolds in a Euclidian space. We actually obtain a sufficient condition for a harmonic foliation on compact submanifolds immersed in Euclidean space to be unstable, where its application to the standard sphere allows us to obtain the result of Kamber and Tondeur \cite{cite8}. We also prove that, for a compact homogeneous Riemannian manifold $(N,g_N)$ satisfying $\lambda_1<2s \cdot \dim N$, any harmonic foliation on $N$ with bundle-like $g_N$ is unstable, where $\lambda_1$ and $s$ denote the first eigenvalue of the Laplacian and the scalar curvature of $g_N$, respectively. Then instability for harmonic foliations on $N$ is equivalent to the non-existence of stable harmonic map between $N$ and any compact Riemannian manifold (see Theorem 4.8). In particular, we determine all simply connected compact irreducible symmetric spaces whose harmonic foliation is unstable (see Theorem 4.9).\vspace{4mm}

This paper is organized as follows. In Section 2, we review the theory of harmonic foliations by Kamber and Tondeur. Then Section 3 is devoted to the proof of Main Theorem above for harmonic foliations. Finally in Section 4, we shall show Theorem 4.8 and 4.9 on instability of harmonic foliations. \vspace{1pc}

{\sc Acknowledgement.} The authors wish to thank Professor Yoshihiro Ohnita for useful suggestion and advice.


\section{The Jacobi operator and a stability of harmonic foliations}

Let $(N,g_N)$ be an $n$-dimensional compact Riemannian manifold and let ${\mathcal F}$ be a foliation given by an integrable subbundle $L\subset TN$. We define a torsion free connection $\nabla$ on normal bundle $Q=TN/L$ by
\begin{equation*}\tag{2.1}
\begin{cases}
\nabla_XS=\pi[X,Y_S], \quad \text{for }X\in\Gamma(L),\ S\in\Gamma(Q) \text{ and } Y_S=\sigma(S)\in\Gamma(\sigma(Q)),\\
\nabla_XS=\pi(\nabla^N_XY_S), \quad \text{for }X\in\Gamma(\sigma(Q)),\ S\in\Gamma(Q) \text{ and } Y_S=\sigma(S)\in\Gamma(\sigma(Q)),
\end{cases}
\end{equation*}
where $\sigma:Q\to TN$ is a splitting such that $\sigma(Q)$ coincides with the orthogonal complement $L^\perp$ of $L$ in $TN$ with respect to $g_N$. If the normal bundle $Q$ is equipped with a holonomy invariant fiber metric $g_Q$, i.e. $Xg_Q(S,T)=g_Q(\nabla_XS,T)+g_Q(S,\nabla_XT)$ for all $X\in\Gamma(TN)$, the foliation $\mathcal F$ is called a {\it Riemannian foliation} or an {\it $R$-foliation}. There is a unique metric $g_Q$ for an $R$-foliation with a torsion free connection $\nabla$ on the normal bundle $Q$. A Riemannian metric $g_N$ on $N$ is called a {\it bundle-like} metric with respect to the foliation $\mathcal F$ if the foliation becomes an $R$-foliation in terms of the fiber metric $g_Q$ induced on $Q$.

For a foliation $\mathcal F$ on a Riemannian manifold $(N,g_N)$, the curvature $R^\nabla$ of the connection $\nabla$ is an $\text{End}(Q)$-valued $2$-form on $N$. Since $i(X)R^\nabla=0$ for $X\in\Gamma(L)$, it follows that the Ricci operator $R^\nabla(S,T):Q\to Q$ for $S,T\in\Gamma(Q)$, is well-defined. Define $P^\nabla(U,V):Q\to Q$ by $P^\nabla(U,V)S=-R^\nabla(U,S)V$ for all $S\in\Gamma(Q)$. The Ricci curvature $S^\nabla$ for $\mathcal F$ is then $S^\nabla(U,V)=\text{trace}P^\nabla(U,V)$ which is a symmetric bilinear form. We define the {\it Ricci operator} $\rho_\nabla:Q\to Q$ as the corresponding self-adjoint operator given by $g_Q(\rho_\nabla U,V)=S^\nabla(U,V)$, where $g_Q$ denotes the holonomy invariant metric on $Q$. In terms of an orthonormal basis $e_{p+1},\ldots,e_n$ of $Q_x$ at some $x\in N$, we have $(\rho_\nabla U)_x=\sum_{\alpha=p+1}^n R^\nabla(U,e_\alpha)e_\alpha$.

Denoting by $\pi\in\Omega^1(N,Q)$ the canonical projection from $TN$ onto $Q$, we have $d_\nabla\pi\in\Omega^2(N,Q)$, $d^*_\nabla\pi\in C^\infty(N, Q)$, the Laplacian $\Delta$ on $\Omega^1(N,Q)$ and so forth. Then we have the following fact (Kamber and Tondeur \cite[3.3]{cite6}).\vspace{1pc}

{\sc Fact.} {\it Let $\mathcal F$ be an $R$-foliation on compact oriented Riemannian manifold $N$ with a bundle-like metric. Then the following are equivalent:}\\
(i) $\pi$ {\it is harmonic}, \\
(ii) {\it all leaves for the foliation are minimal submanifolds of} $N$, \\
(iii) $\Delta\pi=0$.\vspace{1pc}

A foliation is said to be {\it harmonic} if it satisfies (i) or (ii) in above fact.

We next study first and second variations of $R$-foliation ${\mathcal F}$ on a compact Riemannian manifold $(N,g_N)$ with bundle-like metric $g_N$. We define the {\it energy} of the foliation $\mathcal F$ by
$$
E({\mathcal F})=\frac12 \Vert \pi\Vert,
$$
where $\pi$ is the canonical projection from $TN$ onto $Q$ and is considered as a $Q$-valued $1$-form on $N$. Let $\{U_{\alpha},f^\alpha,\gamma^{\alpha\beta}\}$ be the Haefliger cocycle representing ${\mathcal F}$. Namely, $\{U_\alpha\}$ is an open cover of $N$ with $f^\alpha:U_\alpha\to{\mathbb R}^q$ such that $\gamma^{\alpha\beta}$ are local isometries on $U_\alpha\cap U_\beta(\neq\phi)$ satisfying $f^\alpha=\gamma^{\alpha\beta}f^\beta$. Here $q$ denotes the codimension of ${\mathcal F}$. For $\nu\in\Gamma(Q)$, we put
$$
\Phi^\alpha_t(x)=\exp_{f^\alpha(x)}(t\nu^\alpha(x)), \qquad x\in U_\alpha, \ t\in(-\varepsilon,\varepsilon),
$$
where $\nu^\alpha=\nu\vert_{U_\alpha}$. We then have a variation $\Phi^\alpha_t$ of $f^\alpha=\Phi^\alpha_0$, where $\varepsilon$ is sufficiently small. Since $\Phi^\alpha_t(x)=\gamma^{\alpha\beta}\Phi^\beta_t(x)$ on $U_\alpha\cap U_\beta$, the local variations $\{\Phi^\alpha_t\}$ define a variation ${\mathcal F}_t$ of the foliation $\mathcal F$. Moreover we have
\begin{equation*}\tag{2.2}
\nabla_{\frac{\partial}{\partial t}\vert_{t=0}}(\Phi^\alpha_t)_*=\nabla\nu^\alpha\in\Omega^1(U_\alpha,Q).
\end{equation*}
To obtain the second variation, we need a 2-parameter variation ${\mathcal F}_{s,t}$ of ${\mathcal F}_{0,0}={\mathcal F}$ defined locally as $\Phi^\alpha_{s,t}$, where
$$
\Phi^\alpha_{s,t}(x)=\exp_{f^\alpha(x)}(s\mu^\alpha(x))+t\nu^\alpha(x)), \qquad x\in U_\alpha, \ s,t\in(-\varepsilon,\varepsilon)
$$
for $\nu,\mu\in\Gamma(Q)$. Then by (2.2) 
$$
\begin{cases}
\nabla_{\frac{\partial}{\partial s}\vert_{s=0,t=0}}(\Phi_{s,t}^\alpha)_*=\nabla\mu^\alpha, \\
\nabla_{\frac{\partial}{\partial t}\vert_{s=0,t=0}}(\Phi_{s,t}^\alpha)_*=\nabla\nu^\alpha.
\end{cases}
$$
The second variation formula is now given by 
\begin{align*}
& \left.\frac{\partial^2}{\partial s\partial t}\right\vert_{s=0,t=0}E({\mathcal F}_{s,t}) =
\left.\frac{\partial^2}{\partial s\partial t}\right\vert_{s=0,t=0}\frac12 \langle\pi_{s,t},\pi_{s,t}\rangle = \left.\frac{\partial}{\partial s}\right\vert_{s=0,t=0}
\langle\nabla\nu,\pi_{s,t}\rangle \\
& = \langle\nabla_{\frac{\partial}{\partial s}}\nabla\nu,\pi\rangle+\langle\nabla\nu,\nabla\mu\rangle = \langle R^\nabla(\mu,\pi)\nu,\pi\rangle+\langle\nabla\nabla_{\frac{\partial}{\partial s}}\nu,\pi\rangle+\langle d_\nabla\nu,d_\nabla\mu\rangle \\
& = -\langle R^\nabla(\mu,\pi)\pi,\nu\rangle+\langle\nabla_{\frac{\partial}{\partial s}}\nu,d_\nabla^*\pi\rangle+\langle d_\nabla^*d_\nabla\mu,\nu\rangle = \langle(\Delta-\rho_\nabla)\nu,\mu\rangle+\langle\nabla_{\frac{\partial}{\partial s}}\nu,d_\nabla^*\pi\rangle,
\end{align*}
where $R^\nabla$ and $\rho_\nabla$ are the curvature and the Ricci operator for $Q$, respectively. For a harmonic foliation $\mathcal F$, we have
\begin{equation*}\tag{2.3}
\left.\frac{\partial^2}{\partial s\partial t}\right\vert_{s=0,t=0}E({\mathcal F}_{s,t})=\langle(\Delta-\rho_\nabla)\mu,\nu\rangle=\langle{\mathcal J}_\nabla\mu,\nu\rangle,
\end{equation*}
where ${\mathcal J}_\nabla=\Delta-\rho_\nabla$ is the Jacobi operator of $\mathcal F$. Note that the Jacobi operator ${\mathcal J}_\nabla$ is a self-adjoint and strongly elliptic with real eigenvalues $\lambda_1<\lambda_2<\cdots<\lambda_i<\cdots\to\infty$ for $i\to\infty$. Here the dimension of each eigenspace $V_\lambda({\mathcal F})=\{\nu\in\Gamma(Q); {\mathcal J}_\nabla\nu=\lambda\nu\}$ is finite, i.e. $\dim V_\lambda({\mathcal F})<\infty$.\vspace{4mm}

{\sc Definition.} The {\it index} of a harmonic foliatrion $\mathcal F$ is defined by
$$
\text{index}({\mathcal F})=\sum_{\lambda_i<0}\dim V_{\lambda_i}({\mathcal F})
$$
and a harmonic foliation $\mathcal F$ is said to be {\it stable} if $\text{index}({\mathcal F})=0$, i.e. $\langle{\mathcal J}_\nabla\nu,\nu\rangle\geqq0$ for all $\nu\in\Gamma(Q)$.\vspace{4mm}

Note that this definition makes sense for the case of harmonic foliation $\mathcal F$ with bundle-like metric $g_N$, because if $g_N$ is not bundle-like, then the equality (2.3) does not hold in general.


\section{Harmonic foliations on locally conformal K\"ahler manifolds}

The purpose of this section is to prove Main Theorem in Introduction.

For a locally conformal K\"ahler manifold $(M,J,g_M)$ with $\Omega$ and $\omega$, let $B=\omega^\sharp$ be the {\it Lee vector field}, where $\sharp$ denotes the raising of indices with respect ot $g_M$.

The case when $\omega$ is identically zero, $(M,J,g_M)$ is a K\"ahler manifold. Any complex submanifold of a K\"ahler manifold is also K\"ahler, and especially, is minimal. Hence, in this case, we have the following:\vspace{1pc}

\begin{corollary}
The foliations on compact K\"ahler manifolds with a bundle-like metric foliated by complex submanifolds are stable.
\end{corollary}

The following lemma is crucial in the proof of Main Theorem:

\begin{lemma}\label{Lemma32}
The connection $\nabla$ on $Q$ defined in (2.1) satisfies $\nabla_XJ_QS=J_Q\nabla_XS$ for all $X\in\Gamma(TM)$ and $S\in\Gamma(Q)$,
where $J_Q$ denotes the almost complex structure on $Q$ induced by $J$ on $M$.
\end{lemma}

\begin{proof}
We first note that any complex submanifold $N$ of a locally conformal K\"ahler manifold $M$ is minimal if and only if the Lee vector field $B$ for $M$ is tangent to $N$ (for instance, see Dragomir and Ornea \cite[Theorem 12.1]{cite2}).  Let $\nabla^M$ be the Levi-Civita connection of $(M,g_M)$. Then for all $X,Y\in\Gamma(TM)$,
$$
\nabla^M_X JY=J\nabla^M_X Y+\frac12 \{\theta(Y)X-\omega(Y)JX-g_M(X,Y)A-\Omega(X,Y)B\},
$$
where $\theta=\omega\circ J$ and $A=-JB$. Then if $X\in\Gamma(\sigma(Q))$ and $Y\in\Gamma(Q)$, we have
\begin{align*}
\nabla_XJ_QS-J_Q\nabla_XS & = \pi(\nabla^M_XJY_S-J\nabla^M_XY_S) \\
& = \pi(\frac12\{\theta(Y_S)X-\omega(Y_S)JX-g_M(X,Y_S)A-\Omega(X,Y_S)B\}) =0
\end{align*}
by $\theta(Y_S)=\omega(Y_S)=0$. On the other hand, if $X\in\Gamma(L)$ and $S\in\Gamma(Q)$, by Proposition 2.2 of Dragomir and Ornea \cite{cite2} (cf. Vaisman \cite{cite15}), we have $[X,JY_S]-J[X,Y_S]\in L$. Then 
$$
\nabla_XJ_QS-J_Q\nabla_XS=\pi([X,JY_S]-J[X,Y_S])=0,
$$
and this completes the proof of the lemma.
\end{proof}

We define a linear differential operator $D:\Gamma(Q)\to\Gamma(Q\otimes T^*M)$ of first order by
$$
DV(X)=\nabla_{JX}V-J_Q\nabla_XV,\qquad  V\in\Gamma(Q) \text{ and } X\in\Gamma(TM).
$$
{\sc Proof of Main Theorem.} It suffices to show
\begin{equation}\tag{3.3}
\langle{\mathcal J}_\nabla V,V\rangle=\frac12 \langle DV,DV\rangle
\end{equation}
for all $V\in\Gamma(Q)$. Let $\{e_1,\ldots,e_n,f_1,\ldots,f_n\}$ be a local orthonormal frame such that $Je_i=f_i,\ Jf_i=-e_i,\ 1\leqq i\leqq n,$ and that the frame $\{e_1,\ldots,e_p,f_1,\ldots,f_p\}$ spans $\mathcal F$. Then
\begin{align*}\tag{3.4}
\langle {\mathcal J}_\nabla V,V \rangle = & \langle d_\nabla^*d_\nabla V,V \rangle-\langle \rho_\nabla V,V \rangle = \langle d_\nabla V,d_\nabla V \rangle-\langle R^\nabla(V,\pi)\pi,V \rangle \\
= & \sum_{i=1}^n \left\{ \int_M g_Q (\nabla_{e_i}V,\nabla_{e_i}V)v_M+\int_M g_Q(\nabla_{f_i}V,\nabla_{f_i}V) v_M \right\} \\
& -\sum_{i=p+1}^n\left\{ \int_Mg_Q(R^\nabla(V,e_i)e_i,V)v_M+\int_Mg_Q(R^\nabla(V,f_i)f_i,V)v_M \right\}.
\end{align*}
On the other hand, $\langle DV,DV\rangle$ is written as
\begin{align*}\tag{3.5}
\langle DV,DV \rangle = & \sum_{i=1}^n \left\{\int_M g_Q(DV(e_i),DV(e_i))v_M+\int_M g_Q(DV(f_i),DV(f_i))v_M\right\} \\
= & \sum_{i=1}^n \{\int_M g_Q(\nabla_{Je_i}V-J\nabla_{e_i}V,\nabla_{Je_i}V-J\nabla_{e_i}V) \\
& {} \qquad \quad +g_Q(\nabla_{Jf_i}V-J\nabla_{f_i}V,\nabla_{Jf_i}V-J\nabla_{f_i}V)v_M\} \\
= & \sum_{i=1}^n\int_M\{g_Q(\nabla_{Je_i}V,\nabla_{Je_i}V)-2g_Q(\nabla_{Je_i}V,J\nabla_{e_i}V) \\
& {} \qquad \quad +g_Q(J\nabla_{e_i}V,J\nabla_{e_i}V)+g_Q(\nabla_{e_i}V,\nabla_{e_i}V) \\
& {} \qquad \qquad +2g_Q(\nabla_{e_i}V,J\nabla_{Je_i}V)+g_Q(J\nabla_{Je_i}V,J\nabla_{Je_i}V)\}v_M \\
= & 2\sum_{i=1}^n\int_M\{g_Q(\nabla_{e_i}V,\nabla_{e_i}V)+g_Q(\nabla_{Je_i}V,\nabla_{Je_i}V) \\
& {}\qquad \quad + g_Q(\nabla_{e_i}V,J\nabla_{Je_i}V)-g_Q(\nabla_{Je_i}V,J\nabla_{e_i}V)\}v_M \\
= & 2\sum_{i=1}^n\int_M\{g_Q(\nabla_{e_i}V,\nabla_{e_i}V)+g_Q(\nabla_{Je_i}V,\nabla_{Je_i}V)+e_ig_Q(V,J\nabla_{Je_i}V) \\
& -g_Q(V,J\nabla_{e_i}\nabla_{Je_i}V)-Je_ig_Q(V,J\nabla_{e_i}V)+g_Q(V,J\nabla_{Je_i}\nabla_{e_i}V)\}v_M \\
= & 2\sum_{i=1}^n\int_M\{g_Q(\nabla_{e_i}V,\nabla_{e_i}V)+g_Q(\nabla_{Je_i}V,\nabla_{Je_i}V)+ e_ig_Q(V,J\nabla_{Je_i}V) \\
& -Je_ig_Q(V,J\nabla_{e_i}V) - g_Q(V,JR^\nabla(e_i,Je_i)V)-g_Q(V,J\nabla_{[e_i,Je_i]}V)\}v_M.
\end{align*}
We also observe that
\begin{equation*}\tag{3.6}
\sum_{i=1}^n \int_M\{e_ig_Q(V,J\nabla_{Je_i}V)-Je_ig_Q(V,J\nabla_{e_i}V)-g_Q(V,J\nabla_{[e_i,Je_i]}V)\}v_M=0,
\end{equation*}
because if $X\in\Gamma(TM)$ is defined by $g_M(X,Y)=g_Q(\nabla_{JY}V,JV)$, then the following computation of $\text{div}(X)$ together with $\int_M\text{div}(X)v_M=0$ allows us to obtain (3.6):
\begin{align*}
\text{div}(X)= & \sum_{i=1}^n\{g_M(e_i,\nabla_{e_i}^MX)+g_M(Je_i,\nabla_{Je_i}^MX)\} \\
= & \sum_{i=1}^n\{e_ig_M(e_i,X)-g_M(\nabla^M_{e_i}e_i,X)+ Je_ig_M(Je_i,X)-g_M(\nabla^M_{Je_i}Je_i,X)\} \\
= & \sum_{i=1}^n\{e_ig_Q(\nabla_{Je_i}V,JV)-g_Q(\nabla_{J\nabla^M_{e_i}e_i}V,JV) \\
& {}\qquad + Je_ig_Q(\nabla_{JJe_i}V,JV)-g_Q(\nabla_{J\nabla^M_{Je_i}Je_i}V,JV)\} \\
\end{align*}
\begin{align*}
= & \sum_{i=1}^n\{e_ig_Q(\nabla_{Je_i}V,JV)-Je_ig_Q(\nabla_{e_i}V,JV)v\\
& {}\qquad -g_Q(\nabla_{\nabla^M_{e_i}Je_i}V,JV)+g_Q(\nabla_{\nabla^M_{Je_i}e_i}V,JV)\} \\
= & -\sum_{i=1}^n\{e_ig_Q(V,J\nabla_{Je_i}V)-Je_ig_Q(V,J\nabla_{e_i}V)-g_Q(V,J\nabla_{[e_i,Je_i]}V)\}.
\end{align*}
Now by (3.5) and (3.6), we have
\begin{align*}\tag{3.7}
 \langle DV,DV\rangle= & 2\sum_{i=1}^n \int_M\{g_Q(\nabla_{e_i}V,\nabla_{e_i}V)+g_Q(\nabla_{Je_i}V,\nabla_{Je_i}V)\\
& -g_Q(V,JR^\nabla(e_i,Je_i)V)\}v_M.
\end{align*}
Then for $1\leqq i\leqq p$,
\begin{align*}\tag{3.8}
R^\nabla(e_i,Je_i)V & =\nabla_{e_i}\nabla_{Je_i}V-\nabla_{Je_i}\nabla_{e_i}V-\nabla_{[e_i,Je_i]}V \\
& = \pi[e_i,\pi[Je_i,V]]-\pi[Je_i,\pi[e_i,V]]-\pi[[e_i,Je_i],V] \\
& = \pi[e_i,\pi[Je_i,V]]+\pi[Je_i,\pi[V,e_i]]+\pi[V,[e_i,Je_i]]=0,
\end{align*}
because the foliation is involutive satisfying
$$
\pi[e_i,\pi^\perp[Je_i,V]]=0=\pi[Je_i,\pi^\perp[e_i,V]], 
$$
where $\pi^\perp=\text{id}-\pi$. Furthermore, for $p+1\leqq i\leqq n$, the Bianchi identity shows that
\begin{equation*}\tag{3.9}
JR^\nabla(e_i,Je_i)V =-JR^\nabla(Je_i,V)e_i-JR^\nabla(V,e_i)Je_i=R^\nabla(V,Je_i)Je_i+R^\nabla(V,e_i)e_i
\end{equation*}
Thus, by (3.4), (3.7), (3.8) and (3.9), we obtain the required identity (3.3) as follows:
\begin{align*}
\frac12\langle DV,DV\rangle = & \sum_{i=1}^n\int_M\{g_Q(\nabla_{e_i}V,\nabla_{e_i}V)+g_Q(\nabla_{Je_i}V,\nabla_{Je_i}V)\}v_m \\
& - \sum_{i=p+1}^n\int_M\{g_Q(R^\nabla(V,e_i)e_i,V)+g_Q(R^\nabla(V,Je_i)Je_i,V)\}v_M \\
= & \langle{\mathcal J}_\nabla V,V\rangle.
\end{align*}
{}\hfill $\square$\vspace{4mm}

In the remainder of Section 3, we give an example of a stable harmonic foliation on a locally conformal K\"ahler manifold. Let $\lambda$ be a complex number satisfying $\vert \lambda\vert\neq1$. Denote by $\langle\lambda\rangle$ the cyclic group generated by the transformation $:(z_1,\ldots,z_n)\mapsto(\lambda z_1,\ldots,\lambda z_n)$ of ${\mathbb C}^n-\{0\}$. Since this group acts freely and holomorphically on ${\mathbb C}^n-\{0\}$, the quotient space ${\mathbb C}H^n:=({\mathbb C}^n-\{0\})/\langle\lambda\rangle$ is a complex manifold called a {\it Hopf manifold}. Consider the Hermitian metric $g_0=(\Sigma^{n}_{k=1} dz^k\otimes d\bar{z}^k)/\Vert z\Vert^2$ on ${\mathbb C}^n-\{0\}$. Then $g_0$ gives not only a locally conformal K\"ahler structure but also a Vaisman manifold structure on ${\mathbb C}H^n$ with Lee form $\omega_0=-\{\Sigma^{n}_{k=1}(z^kd\bar{z}^k+\bar{z}^kdz^k)\}/\Vert z\Vert^2$. It is well-known that ${\mathbb C}H^n$ has a principal $T^1_{\mathbb C}$-bundle structure over the projective space ${\mathbb C}P^{n-1}$. Then the foliation on ${\mathbb C}H^n$ defined by the canonical projection $\pi:{\mathbb C}H^n\to{\mathbb C}P^{n-1}$ is harmonic and is stable by Main Theorem, where the metric on ${\mathbb C}P^{n-1}$ is the Fubini-Study metric. 

\begin{remark0}
(i) More generally, Main Theorem is valid even if $M$ is (not necessarily K\"ahler and is) just a compact Hermitian manifold, provided that the connection $\nabla$ defined by (2.1) satisfies Lemma 3.2.\\
(ii) As to stable harmonic foliations, there exists an example foliated by fibers of a Riemannian submersion whose base space is not a complex manifold . A typical example is the twistor space of a quaternionic K\"ahler manifold.
\end{remark0}


\section{Instability of harmonic foliations}

In this section, we discuss instability for harmonic foliations on Riemannian manifolds. Let $(N,g_N)$ be a Riemannian manifold. By the Weitzenb\"ock formula we have $(\Delta\pi)=\nabla^*\nabla\pi+S(\pi)$, where by using a local frame $\{e_1,\ldots,e_n\}$ for $TN$, we put
$$
\nabla^*\nabla\pi=-\sum_{i=1}^n(\nabla^2_{e_i,e_i}\pi) \text{ and } S(\pi)(X)=\sum_{i=1}^n \{R^\nabla(e_i,X)\pi(e_i)-\pi(R^N(e_i,X)e_i)\}
$$
for all $X\in\Gamma(TN)$. Here $R^\nabla$ and $R^N$ denote the curvature tensors associated to $\nabla$ and $\nabla^N$, respectively. We then have
\begin{equation*}\tag{4.1}
\Delta\pi=\nabla^*\nabla\pi-\rho_\Delta\cdot\pi+\pi\cdot\rho_N,
\end{equation*}
in view of the equality
$$
S(\pi)(X) = -\sum_{i=p+1}^n R^\nabla(\pi(X),e_i)e_i+\pi(\sum_{i=1}^n R^N(X,e_i)e_i)= -(\rho_\Delta\pi(X))+\pi(\rho_N(X)).
$$

Let $\mathcal F$ be a Riemannian and harmonic foliation on $N$ with bundle-like $g_N$, i.e., the canonical projection $\pi:TN\to Q$ satisfies $\Delta\pi=0$. Then (4.1) is expressible as
\begin{equation*}\tag{4.2}
\rho_\nabla\cdot\pi=\nabla^*\nabla\pi+\pi\cdot\rho_N.
\end{equation*}
On the other hand, by operating the Laplacian on $\pi(X), X\in\Gamma(TN),$ we obtain
\begin{align*}\tag{4.3}
\Delta(\pi(X)) & = d^*_\nabla d_\nabla(\pi(X))=-\sum_{i=1}^n \nabla^2_{e_i,e_i}(\pi(X)) \\
& = (\nabla^*\nabla)(X)+\pi(\nabla^{N*}\nabla^NX)-2\sum_{i=1}^n(\nabla_{e_i}\pi)(\nabla^N_{e_i}X).
\end{align*}
By (4.2) and (4.3), ${\mathcal J}_\nabla(\pi(X))=(\Delta-\rho_\nabla)(\pi(X))$ is written as
$$
{\mathcal J}_\nabla(\pi(X))=-(\pi\cdot\rho_N)(X)+\pi(\nabla^{N*}\nabla^NX)-2\sum_{i=1}^n (\nabla_{e_i}\pi)(\nabla^N_{e_i}X).
$$
Assume that $N$ is a compact submanifold immersed in the Euclidean space ${\mathbb E}^n$ with the standard inner product $\ll\ , \gg$. For each vector $v$ in ${\mathbb E}^n$, we define a smooth function $f_v$ on $N$ by $f_v(x):=\ \ll v,x\gg$ for $x\in N$. We denote by $\Psi_t, t\in{\mathbb R},$ the flow generated by $V=\text{grad}f_v$. Simple computations give us
\begin{equation*}\tag{4.4}
\ll \nabla^N_X,Y \gg \ =\ \ll B(X,Y),v\gg,
\end{equation*}
\begin{equation*}\tag{4.5}
\ll (\nabla^N)^2_{X,Z},Y\gg\ =\ -\ll B(X,Y),B(Z,V)\gg+\ll(\nabla B)(X,Y,Z),v\gg,
\end{equation*}
where $B$ denotes the second fundamental form for the submanifold $N$ in ${\mathbb E}^n$.

The energy functional for $\mathcal F$ is defined by $E({\mathcal F})=(1/2) \int_N\Vert\pi\Vert$. Consider the associated quadratic form $Q_{\mathcal F}$ by setting
$$
Q_{\mathcal F}(v)=\frac{d^2}{dt^2}E(\Psi_t)\vert_{t=0}=\int_N g_N({\mathcal J}_\nabla(\pi(V)),\pi(V)).
$$
We shall now compute the trace $\text{Tr}(Q_{\mathcal F})$ of $Q_{\mathcal F}$ on ${\mathbb E}^n$. By (4.4) and (4.5),
\begin{align*}
g_Q ({\mathcal J}_\nabla(\pi(V)),\pi(V))= & -g_Q(\pi(\rho_N(V)),\pi(V)) \\
 & +\sum_{k,l=1}^n \ll B(e_k,e_l),B(e_k,V)\gg g_Q(\pi(e_l),\pi(V)) \\
& -\sum_{k,l=1}^n \ll (\nabla B)(e_k,e_l,e_k),v\gg g_Q(\pi(e_l),\pi(V)) \\
& -2\sum_{k,l=1}^n \ll B(e_k,e_l),v\gg g_Q((\nabla_{e_k}\pi)(e_l),\pi(V)).
\end{align*}
Hence we have
\begin{align*}
\text{Tr}(Q_{\mathcal F}) & = \int_N\{-\sum_{k=1}^n g_Q(\pi(\rho_N(e_k)),\pi(e_k)) \\
& {}\qquad +\sum_{k,l,m=1}^n \ll B(e_k,e_l),B(e_l,e_m)\gg g_Q(\pi(e_k),\pi(e_m))\} \\
& =\int_N\sum_{a=p+1}^n\{\sum_{j=1}^p \ll B(e_a,e_j),B(e_j,e_a)\gg-\ll \rho_N(e_a),e_a\gg\}.
\end{align*}
Let $\eta$ denote  the mean curvature vector of the submanifold $N$ in ${\mathbb E}^n$. Then by the equation of Gauss, we obtain
\begin{equation*}
\text{Tr}Q_{\mathcal F}=\int_N\sum_{a=p+1}^n(n\ll B(e_a,e_a),\eta\gg-2\ll\rho_N(e_a),e_a\gg).
\end{equation*}
This immediately implies \vspace{4mm}

{\sc Lemma 4.6.} {\it  Let $(N,g_N)$ be an $n$-dimensional compact submanifold immersed in the Euclidean space ${\mathbb E}^N$. If $N$ satisfies
$$
n\ll B(u,u),\eta\gg-2\ll \rho_N(u),u\gg\ <0
$$
for all unit vector $u$ in $TN$, then every Riemannian and harmonic foliation on $N$ with bundle-like $g_N$ is unstable.}\vspace{4mm}

In the case where $N$ is the standard sphere, the above result was proved by Kamber and Tondeur \cite{cite8}. This lemma is a generalization of a result of Ohnita \cite{cite10} known for harmonic maps. The following is now straightforward from Lemma 4.6. \vspace{4mm}

{\sc Theorem 4.7.} {\it  Let $N$ be an $n$-dimensional compact minimal submanifold of a unit sphere $S^{N-1}(1)$. If the Ricci curvature $S_N$ of $N$ satisfies $S_N>2/n$, then every Riemannian and harmonic foliation on $N$ with bundle-like $g_N$ is unstable.}\vspace{4mm}

It might be of some interest to compare results on instability for harmonic foliations with that of harmonic maps. Hence, by combining Theorem 4.7 above with Theorem 4 of \cite{cite10}, we obtain: \vspace{4mm}

{\sc Theorem 4.8.} {\it Let $(N,g_N)$ be an $n$-dimensional compact homogeneous Riemannian manifold with irreducible isotropy representation. For $(N,g_N)$, let $s$ and $\lambda_1$ denote the scalar curvature and the first eigenvalue of the Laplacian acting on functions, respectively. Then the following conditions are all equivalent:\\
(1) $\lambda_1<2s/n$. \\
(2) Every Reimannian and harmonic foliation on $N$ with bundle-like $g_N$ is unstable. \\
(3) There exist no nonconstant stable harmonic maps from $N$ to Riemannian manifolds. \\
(4) There exist no nonconstant stable harmonic maps from compact Riemannian manifolds to $N$. \\
(5) The identity map $\text{id}_N$ of $N$ onto itself is unstable as a harmonic map.}\vspace{4mm}

\begin{proof}
The implications (3) $\Longrightarrow$ (5) and (4) $\Longrightarrow$ (5) are trivial. Since the stability of the point foliation on $N$ is equivalent to the stability of $\text{id}_N$ as a harmonic map, (2) implies (5). Since $N$ is an Einstein manifold from a result of Smith \cite{cite11}, we have the equivalence (1) $\Longleftrightarrow$ (5). Hence, it suffices to show (1) implies (2), (3) and (4). By virtue of the theorem of Takahashi \cite{cite12}, there exists a standard minimal immersion $\varphi$ of $N$ into a unit hypersphere $S^m(1)$ by using an orthonormal basis for the first eigenspace of the Laplacian in such a way that $\varphi$ is an isometric immersion of $(N,(\lambda_1/n)g_N)$ into $S^m(1)$. Then the Ricci curvature of $(N,(\lambda_1/n)g_N)$ is greater than $n/2$. By Theorem 4.7 and Theorem 1 of Ohnita \cite{cite10}, we obtain (2), (3) and (4).
\end{proof}

\begin{remark0}
Theorem 4.8 is valid even if we replace homogeneous $N$ above by a strongly harmonic manifold. However, for strongly harmonic manifolds, no inhomogeneous examples are known (c.f. Besse \cite{cite1}).
\end{remark0}

Compact irreducible symmetric spaces which satisfy $\lambda_1<2s/n$ were determined by Smith \cite{cite11}, Nagano \cite{cite9} and Ohnita \cite{cite10}. Thus we obtain\vspace{4mm}

{\sc Theorem 4.9.} {\it Let $(N,g_N)$ be a compact irreducible symmetric space. Then the following conditions are equivalent:\\
(1) Any Riemannian and harmonic foliation on $N$ with bunde-like $g_N$ is unstable. \\
(2) $N$ is simply connected and belongs to one of the following:\\
 (a) $SU(n)\ (n\geqq2)$, \quad (b) $Sp(n) \ (n\geqq2)$ \quad (c) $SU(2n)/Sp(n) \ (n\geqq3)$ \quad (d) $S^n \ (n\geqq 3)$ \\
 (e) $G_{p,q}({\mathbb H})=Sp(p+q)/Sp(p)\times Sp(q) \ (p\geqq q \geqq 1)$ \quad (f) $E_6/F_4$ \quad (g) $P_2({\mathbb O})=F_4/Spin(9)$.}\vspace{4mm}

Applying Lemma 4.6 to product isometric immersion, we have the following:\vspace{4mm}

{\sc Corollary.} {\it If $(N,g_N)$ is a product of simply connected compact irreducible symmetric spaces belonging to the list in (2) of Theorem 4.9, then every Riemannian and harmonic foiation on $N$ with bundeli-like $g_N$ is unstable.}



\begin{thebibliography}{99}
\bibitem{cite1}
\textsc{A. Besse},
Manifolds all of whose geodesics are closed,
 With appendices by D. B. A. Epstein, J.-P. Bourguignon, L. Berard-Bergery, M. Berger and J. L. Kazdan. Ergebnisse der Mathematik und ihrer Grenzgebiete [Results in Mathematics and Related Areas], 93. Springer-Verlag, Berlin-New York, 1978.

\bibitem{cite2}
\textsc{S. Dragomir and L. Ornea},
Locally conformal K\"ahler geometry,
Progress in Math. 155, Birkh\"auser, Boston, 1998.

\bibitem{cite3}
\textsc{J. Eells and L. Lemaire},
Selected topics of Harmonic maps,
CBMS Regional Conference Series in Mathematics, 50. Published for the Conference Board of the Mathematical Sciences, Washington, DC; by the American Mathematical Society, Providence, RI, 1983.

\bibitem{cite4}
\textsc{J. Eells and J.H. Sampson},
Harmonic mapping of Riemannian manifolds, Amer. J. Math. 86 (1964), 109--160.

\bibitem{cite5}
\textsc{S. Kobayashi and K. Nomizu},
Foundation of differential geometry I, II, John Wiley • Sons, Inc., New York.

\bibitem{cite6}
\textsc{F.W. Kamber and Ph. Tondeur},
Harmonic Foliation, Lecture Note in Math. 949, 87--121, Springer-Verlag, Berlin, Heidelberg, New York, 1982.

\bibitem{cite7}
\textsc{F.W. Kamber and Ph. Tondeur},
Infinitesimal automorphisms and second variation of the energy for harmonic foliation, T\^ohoku Math. J., 34 (1982), 525--538.

\bibitem{cite8}
\textsc{F.W. Kamber and Ph. Tondeur},
The index of harmonic foliations on spheres, Trans. Amer. Math. Soc. 275 (1983), 257--263.

\bibitem{cite9}
\textsc{T. Nagano}
Stability of harmonic maps between symmetric spaces, Lecture Notes in Math. 949, 130--137, Springer-Verlag, New York, 1982.

\bibitem{cite10}
\textsc{Y. Ohnita},
Stability of harmonic maps and standerd minimal immersions, T\^ohoku Math. J. 38 (1986), 259--267.

\bibitem{cite11}
\textsc{R.T. Smith},
The second variational formulas for harmonic mappings, Proc. Amer. Math. Soc. 47 (1975), 229--236.

\bibitem{cite12}
\textsc{T. Takahashi},
Minimal immersions of Riemannian manifolds, J. Math. Soc. Japan 18 (1966), 380-385.

\bibitem{cite13}
\textsc{Ph. Tondeur},
Foliation on Riemannian manifolds, Springer-Verlag, New York, 1988

\bibitem{cite14}
\textsc{H. Urakawa},
Calculus of variation and harmonic map, Translated from the 1990 Japanese original by the author. Translations of Mathematical Monographs, 132. American Mathematical Society, Providence, RI, 1993.

\bibitem{cite15}
\textsc{I. Vaisman},
On the analytic distributions and foliations of a Kaehler manifold, Proc. A.M.S. 58 (1976), 221--228.
\end{thebibliography}
\end{document}